%% file: semi.tex
\documentclass[12pt]{amsart}
\usepackage{amssymb, euscript, diagrams}
\makeatletter \@addtoreset{equation}{section} \makeatother

\newcommand{\xref}[1]{{\rm \ref{#1}}}

\def\stackunder#1#2{\mathrel{\mathop{#2}\limits_{#1}}}%
\newcommand{\muu}{\mbox{\boldmath $\mu$}}
\newcommand{\Diff}{\operatorname{Diff}}
\newcommand{\Pic}{\operatorname{Pic}}
\newcommand{\qq}{\mathbin{\sim_{\scriptscriptstyle{\QQ}}}}
\renewcommand{\emptyset}{\varnothing}

\newcommand{\CC}{\mathbb{C}}
\newcommand{\ZZ}{\mathbb{Z}}
\newcommand{\PP}{\mathbb{P}}
\newcommand{\NN}{\mathbb{N}}
\newcommand{\QQ}{\mathbb{Q}}
\newcommand{\OOO}{\mathcal{O}}
\newcommand{\RRR}{\mathcal{R}}
\newcommand{\De}{\Delta}
\newcommand{\Ga}{\Gamma}
\newcommand{\var}{\varphi}
\newcommand{\ep}{\varepsilon}
\newcommand{\Supp}{\operatorname{Supp}}
\newcommand{\Sing}{\operatorname{Sing}}

\newcommand{\red}{\operatorname{red}}
\newcommand{\T}{\operatorname{T}}
\newcommand{\mt}[1]{\operatorname{#1}}

\newcommand{\down}[1]{\left\lfloor #1\right\rfloor}
\newcommand{\up}[1]{\left\lceil #1\right\rceil}

\newtheorem{theorem}[equation]{Theorem}
\newtheorem{proposition}[equation]{Proposition}
\newtheorem{lemma}[equation]{Lemma}
\newtheorem{corollary}[equation]{Corollary}

\theoremstyle{definition}
\newtheorem{example}[equation]{Example}
\newtheorem{definition}[equation]{Definition}
\newtheorem{notation}[equation]{Notation}
\newtheorem{construction}[equation]{Construction}
\newtheorem{pusto}[equation]{}
\newtheorem{remark}[equation]{Remark}

\theoremstyle{remark}

\author{Yuri Prokhorov}
\title{On semistable Mori contractions}
\address{Department of Algebra,
Faculty of Mathematics, Moscow State University, Moscow 117234,
Russia} \email{prokhoro@mech.math.msu.su}

\date{}
\begin{document}

\begin{abstract}
We study Fano-Mori contractions with fibers of dimension at most
one satisfying the semistability assumption. As an application of
our technique we give a new proof of the existence of semistable
$3$-fold flips.
\end{abstract}
\maketitle
\section{Introduction}
This paper is a continuation of our study of Mori contractions
from threefolds to surfaces (see \cite{P-sb}, \cite{P2},
\cite{P-rb}). We refer to \cite{Ut}, \cite{KM2} for the
terminology of the minimal model theory.

Let $X$ be a normal algebraic threefold over $\CC$ (or
three-dimensional normal complex space) with only terminal
singularities. A proper surjective morphism $f\colon X\to Z$ is
called a \emph{Fano-Mori contraction} if $f_*\OOO_X=\OOO_Z$ and
the anticanonical divisor $-K_X$ is $f$-ample. Our interest is in
the \emph{local structure} of such contractions, so we shall
always assume that $Z$ is not a point and $Z$ and $X$ are
sufficiently small (algebraic or analytic) neighborhoods of some
point $o\in Z$ and the fiber $f^{-1}(o)$, respectively. Note that
we do not assume that $f$ is an \emph{extremal Mori contraction}
(i.e., $X$ is $\QQ$-factorial and $\rho(X/Z)=1$). If the dimension
of fibers of $f$ (near $f^{-1}(o)$) is at most one we can
distinguish the following cases:
\par\smallskip
$\scriptstyle\bullet\ $ $\dim Z=2$, then $f$ is called a
\emph{Mori conic bundle},

$\scriptstyle\bullet\ $ $\dim Z=3$ and $f$ contracts a divisor to
a curve, then $f$ is called a \emph{$2$-$1$-type divisorial
contraction},

$\scriptstyle\bullet\ $ $\dim Z=3$ and the exceptional locus of
$f$ is one-dimensional, then $f$ is called a \emph{flipping
contraction}.
\par\smallskip\noindent
In this paper we deal with \emph{semistable} Fano-Mori
contractions which appear in the semistable minimal model program,
see \cite{Tsunoda}, \cite{Kaw}, \cite{Sh1}, \cite{Kaw1},
\cite{KM2}, and references therein.

\begin{definition}
\label{def-semi}
A Fano-Mori contraction $f\colon X\to Z\ni o$ is said to be
\emph{semistable} if there exists an effective Cartier divisor
$o\in T\subset Z$ such that $(X,f^*T)$ is divisorial log terminal
(dlt).
\end{definition}

It is clear that in the above definition we may replace $T$ with a
general hyperplane section through $o$. In particular, in the case
$\dim Z=3$ we may assume that $T$ is irreducible, normal, and
$(X,f^*T)$ is purely log terminal (plt).

\begin{example}[{\cite[Example 2.1]{P-sb}}]
\label{ex-toric}
The toric contraction $(\PP^1\times \CC^2)/\muu_n(0:a;1,-1)\to
\CC^2/\muu_n(1,-1)$ is a semistable Mori conic bundle with
$T=\{x_1x_2=0\}/\muu_n$.
\end{example}

We shall show that the example above is very special: in ``most''
cases the surface $f^{-1}(T)$ is irreducible and normal.

\begin{proposition}
\label{prop-1-2}
Let $f\colon X\to Z\ni o$ be a semistable Mori conic bundle and
let $T$ be a general hyperplane section through $o$.
\begin{enumerate}
\item
If $T$ is reducible, then $f$ is analytically isomorphic to the
Mori conic bundle from Example \xref{ex-toric}. In particular,
$Z\ni o$ is a Du Val point of type $A_{m-1}$.
\item
If $T$ is irreducible, then $Z\ni o$ is smooth and the pair
$(X,f^*T)$ is purely log terminal \textup(plt\textup). In this
case $S:=f^*T$ is a normal surface with cyclic quotient
singularities of type $\T$ or Du Val of type $A_n$ \textup(see \S
\xref{sect-T} for the definition of $\T$-singularities\textup).
\end{enumerate}
\end{proposition}
In case (ii) the structure of $X$ and $f$ is completely determined
by the structure of the surface $S$ and the contraction $S\to T$.
We study such contractions in \S \xref{sect-lcb}.

In practice, it is very difficult to construct nontrivial examples
of Mori conic bundles explicitly (cf. \cite[\S 5]{P-sb}). Using
deformation theory it is very easy to prove the existence (or
non-existence) of semistable ones, see \S \xref{S-def}. In
particular, we prove the following.

\begin{theorem}
For any three-dimensional terminal semistable singularity $U\ni
P$, there is a semistable Mori conic bundle $f\colon X\to Z\ni o$
with a unique singular point which is analytically isomorphic to
$U\ni P$.
\end{theorem}

The following result was inspired by M. Reid's ``general
elephant'' conjecture (cf. \cite[\S 4]{P-sb}) and provides an
evidence for it.

\begin{theorem}[cf. {\cite[Th. 0.4.5]{Mori}} {\cite[Th. 2.2]{KM1}},
{\cite[Corr. 4.9]{Sh1}}]
\label{ge}
Let $f\colon X\to Z\ni o$ be a semistable Fano-Mori contraction
such that the dimension of fibers is at most one and let $T$ be a
general hyperplane section through $o$. Then $K_X+f^*T$ is
$1$-complementary \cite{Lect}, i.e., there exists an effective
integral Weil divisor $F$ such that $K_X+f^*T+F$ is log canonical
\textup(lc\textup) and linearly trivial over $Z$. Moreover,
$K_X+F$ is canonical and linearly trivial, the surface $F$ is
normal, has only Du Val singularities of type $A_n$, and in the
Stein factorization $F\to \bar F\to f(F)$, the same holds for
$\bar F$.
\end{theorem}
Using this theorem we give a new proof of the existence of
semistable flips in \S \xref{sect-flips}. Our proof is based on
Theorem \xref{ge}, Kawamata's double covering trick \cite{Kaw},
and the existence of certain canonical flops.

\subsection*{Terminology}
The semistable MMP is originated in semistable degenerations of
surfaces. Namely, let $h\colon \mathfrak{X}\to \varDelta$ be a
projective surjective morphism from smooth threefold to a smooth
curve such that the general fiber is a smooth surface and special
fibers are reduced simple normal crossings divisors. In this
situation $\mathfrak{X}/\varDelta$ satisfies the following
property:
\begin{quote}
(*)\qquad\qquad $(\mathfrak{X},h^*P)$ is dlt for every point $P\in
\varDelta$.
\end{quote}
In order to obtain either a minimal or relative Fano model we run
the $K$-MMP over $\varDelta$. Every step of the $K$-MMP is in the
same time a step of the $K+h^*P$-MMP. In particular, the property
(*) is preserved and all contractions and flips are semistable in
our sense. This agrees more or less with definitions given in
\cite{Tsunoda}, \cite{Kaw}, \cite{Corti}.

The semistability defined by Shokurov in \cite{Sh1} (cf.
\cite{Kulikov}) is also close to our one by \cite[Lemma 1.4]{Sh1}.
However the construction is inductive and given in terms of some
(not necessarily projective) resolution.

Koll\'ar and Mori \cite[p. 541]{KM1} defined semistable extremal
neighborhoods $f\colon X\to Z\ni o$ in terms of general member
$F_Z\ni |-K_Z|$: $f$ is semistable if $F_Z\ni o$ is a Du Val
singularity of type $A_n$. By Theorem \xref{ge} our definition
\xref{def-semi} implies Mori-Koll\'ar's one. Conversely, if
$F_Z\ni o$ is a Du Val singularity of type $A_n$, then
$K_X+F+f^*T$ is log canonical (but not necessarily dlt) for some
effective Cartier divisor $T$ on $Z$. Thus $f$ is ``almost''
semistable in our sense.

Our technique uses the Kawamata-Viehweg vanishing theorem, so the
proofs work only in characteristic zero. The positive and mixed
characteristic case was treated in \cite{Kaw1}.

\subsection*{Acknowledgments}
I would like to thank Kyoto Research Institute for Mathematical
Sciences for the hospitality during my stay there in 2002--2003.
This paper is based on the subject of my talk given on the
algebraic geometry seminar of RIMS on April 17, 2003. I am very
grateful to the participants of this seminar for their attention
and V. V. Shokurov for useful comments. The work was partially
supported by the grant INTAS-OPEN-2000-269.

\section{Preliminary results}
In this section we prove Proposition \xref{prop-1-2}.

Let $S$ be an algebraic surface (or two-dimensional complex space)
having at worst quotient singularities and let $\var\colon S\to T$
be a contraction. We assume that $T$ is not a point and $T$ is a
sufficiently small neighborhood of $o\in T$. We say that $\var$ is
a \emph{log contraction} if $-K_S$ is $\var$-ample. If furthermore
any singularity of $S$ is of type $\T$ or Du Val, then we say that
$\var$ is a \emph{$\T$-contraction}.

\begin{lemma}[see {\cite[Lemma 2.5]{P-rb}}]
\label{lemma-plt-plt}
Let $f\colon X\to Z\ni o$ be a Mori conic bundle and let $T$ be an
effective $\QQ$-Weil divisor on $Z$ such that $(X,f^*T)$ is lc
\textup(resp. plt\textup) at some point $P\in f^{-1}(o)$. Then
$(Z,T)$ is lc \textup(resp. plt\textup).
\end{lemma}

\begin{lemma}
Let $f\colon X\to Z\ni o$ be a semistable Mori conic bundle and
let $T$ be an effective Cartier divisor such that $(X,f^*T)$ is
dlt. If $T$ is irreducible, then $Z\ni o$ is smooth and for a
general hyperplane section $o\in T_{\mt{gen}}\subset Z$ the pair
$(X,f^*T_{\mt{gen}})$ is plt.
\end{lemma}
\begin{proof}
Follows by Lemma \xref{lemma-plt-plt} and Bertini's theorem.
\end{proof}

\begin{proposition}
\label{prop-dlt-2}
Let $f\colon X\to Z\ni o$ be a Mori conic bundle. Assume that
there is an effective Weil divisor $T$ on $Z$ such that $(X,f^*T)$
is lc. If $T$ is reducible, then $f$ is analytically isomorphic to
one of the following the Mori conic bundles
\begin{enumerate}
\item
$f$ from Example \xref{ex-toric}, or
\item
$X'/\muu_2(1:0:0;1,1)\to \CC^2/\muu_2(1,1)$, where $X'=\{
x_0^2+x_1^2+x_2^2\phi (u, v) =0\}\subset \PP^2_{x_0, x_1,
x_2}\times\CC^2_{u, v}$ and $\phi (u,v)$ is a $\muu_2$-invariant
without multiple factors, see \cite[Example 2.3]{P-sb}.
\end{enumerate}
In particular, $Z\ni o$ is Du Val of type $A_{m-1}$. Moreover, the
statement of Theorem \xref{ge} holds for $f$.
\end{proposition}
\begin{proof}
By Lemma \xref{lemma-plt-plt}, $T$ has exactly two components
$T_1$ and $T_2$. Consider a base change
\begin{equation*}
\begin{diagram}[height=20pt,width=35pt]
X'&\rTo^{\scriptstyle\upsilon}&X
\\
\dTo^{\scriptstyle {f'}}&&\dTo{\scriptstyle {f}}
\\
Z'&\rTo^{\scriptstyle \pi}&Z
\end{diagram}
\end{equation*}
where $Z'$ is smooth, $Z=Z'/\muu_n$, $X=X'/\muu_n$, and $\muu_n$
acts on $Z'$ free in codimension one (see \cite[(1.9)]{P-sb}). Put
$S=f^*T$, $S'=\upsilon^{-1}(S)$, $T'=\pi^*(T)$, $T'_i=\pi^*(T_i)$,
and $S_i'=f'^*T_i'$. Then $(X',S')$ is lc. Replacing $T_1'$ and
$T_2'$ with general hyperplane sections, we may assume that $T'_1$
is smooth and $(S_1',{S_1'}|_{S_2'})$ is lc, where
${S_1'}|_{S_2'}$ is a Cartier divisor on $S_1'$. In this
situation, $S_1'$ has at worst Du Val singularities. Hence $X'$ is
Gorenstein. Since $f^{\prime -1}(o')$ is reduced, by \cite[\S
2]{P-sb} $X/Z$ we have only two choices for the action of
$\muu_n$. Finally, the statement of Theorem \xref{ge} easily
follows by the proposition below.
\end{proof}

\begin{proposition}[see {\cite[Prop. 4.4.1]{Lect}}]
\label{prop-compl-ext}
Let $f\colon X\to Z\ni o$ be a contraction and let $S$ be a
reduced divisor on $X$ such that $S\cap f^{-1}(o)\neq\emptyset$,
$(X,S)$ is plt and $-(K_X+S)$ is $f$-nef and $f$-big. If
$K_S+\Diff_S$ is $n$-complementary, then so is $K_X+S$. Here
$\Diff_S$ is the \emph{different}, a correcting term in the
Adjunction Formula $K_S+\Diff_S=(K_X+S)|_S$, see \cite[Ch.
16]{Ut}.
\end{proposition}

From now on we consider semistable Fano-Mori contractions $f\colon
X\to Z\ni o$ such that $(X,f^*T)$ is plt (and the dimension of
fibers is at most one).

\section{Singularities of class $\T$}
\label{sect-T}
Let $\muu_n$ acts on $\CC^2$ via $(x,y)\to (\eta^a x,\eta^b y)$,
where $\eta$ is a primitive $n$th root of unity and
$\gcd(n,a)=\gcd(n,b)=1$. In this case we say that the quotient
$\CC^2/\muu_n$ is a singularity of type $\frac1n(a,b)$. This
singularity can be written as $\frac1n(1,q)$, so it is determined
by the fraction $n/q$. The minimal resolution of $\frac1n(1,q)$
can be described as follows. The dual graph of the exceptional
divisor is a chain of smooth rational curves whose
self-intersections $-b_1,\dots,-b_\varrho$ are determined by the
continued fraction expansion
\begin{equation}
\label{eq-cont-frac}
\frac nq=
b_1-\cfrac{1}{b_2-\cfrac{1}{\cdots\cfrac{1}{b_\varrho}}}.
\end{equation}
For typographical reasons we denote the fraction in
\eqref{eq-cont-frac} by $[b_1,\dots,b_\varrho]$. Put
$\varrho(n/q):=\varrho$. Define also the following invariants:
\begin{itemize}
\item
$\iota(n/q)=n/\gcd(n,q+1)$, the index of $\frac1n(1,q)$,
\item
$\beta(n/q)=\gcd(n,q+1)/\iota(n/q)=\gcd(n,q+1)^2/n$,
\item
$\gamma(n/q)=(q+1)/\gcd(n,q+1)$.
\end{itemize}
By definitions, $\iota, \gamma\in \NN$, $\gamma\le \iota$,
$\gcd(\iota, \gamma)=1$. Thus we have the triple
$(\iota,\beta,\gamma)$ which determines $n/q$:
\[
n=\beta\iota^2,\quad q=\beta\iota\gamma-1,\quad
\]
Note that presentation of a cyclic quotient singularity in the
form $\frac1n(1,q)$ is not unique: $\frac1n(1,q')$ defines the
same singularity if and only if either $q\equiv q'$ or $qq'\equiv
1 \mod n$. Clearly, $\varrho(n/q)=\varrho(n/q')$. Since
$\gcd(n,q+1)=\gcd(n,q'+1)$, we have $\beta(n/q)=\beta(n/q')$ and
$\iota(n/q)=\iota(n/q')$.

\begin{definition}
A cyclic quotient non-Du Val singularity $\frac1n(a,b)$ is said to
be of type $\T$ (or simply $\T$-singularity) if $(a+b)^2\equiv 0
\mod n$. (It is easy to see that this definition does not depend
on the representation in the form $\frac1n(a,b)$).
\end{definition}
If $\frac 1n(1,q)$ is a $\T$- (resp. Du Val) singularity, then we
say that $n/q$ is a $\T$- (resp. Du Val) fraction and
$[b_1,\dots,b_\varrho]$ is $\T$- (resp. Du Val) chain. Thus $n/q$
is a $\T$-fraction if and only if $\beta(n/q)\in \ZZ$.

\begin{lemma}
\label{lemma-q-q1}
Let $qq'\equiv 1\mod n$. Then $n/q$ is a $\T$-fraction if and only
if $q+q'=n-2$ if and only if
$\gamma(n/q')+\gamma(n/q)=\iota(n/q)$.
\end{lemma}

\begin{remark}
If $\iota(n/q)=2$, then $\gamma(n/q)=1$ and $n/q$ is a
$\T$-fraction. It is easy to see that this implies either
\begin{equation}
\label{eq-graph-basic}
n/q=[4], \qquad \text{or}\qquad n/q=[3,2,\dots,2,3].
\end{equation}
Moreover, $\beta(n/q)=\varrho(n/q)$. Conversely, any chain such as
in \eqref{eq-graph-basic} has $\iota=2$.
\end{remark}

The minimal resolutions of $\T$-singularities are completely
described.

\begin{proposition}[\cite{KS}]
\label{prop-T-ind}
\begin{enumerate}
\item
If the chain $[b_1,\dots,b_\varrho]$ is of type $\T$, then so are
the chains
\begin{equation*}
a)\quad [2,b_1,\dots,b_\varrho+1] \quad\text{and}\quad b)\quad
[b_1+1,\dots,b_\varrho,2]
\end{equation*}
\item
Every $\T$-chain can be obtained by starting with one of the
chains \eqref{eq-graph-basic} and iterating the steps described in
{\rm (i)}.
\end{enumerate}
\end{proposition}

For a chain $[b_1,\dots,b_\varrho]$, we denote corresponding log
discrepancies by $\alpha_1,\dots,\alpha_\varrho$.

\begin{lemma}
In the above notation one has $\alpha_1= (q+1)/n= \gamma/\iota$.
If $n/q$ is a $\T$-fraction, then $\alpha_\varrho=1-\gamma/\iota$.
\end{lemma}
\begin{proof}
The $\frac1n(1,q)$-weighted blow-up gives us the first relation.
The second one follows by Lemma \xref{lemma-q-q1}.
\end{proof}

\begin{corollary}
\label{cor-alpha-alpha}
Let $[b_1,\dots,b_\varrho]$ be any chain. The following are
equivalent:
\begin{enumerate}
\item
$[b_1,\dots,b_\varrho]$ is of type $\T$;
\item
$\alpha_1+\alpha_\varrho=1$.
\end{enumerate}
\end{corollary}

\begin{theorem}[{\cite[Prop. 5.9]{LW}}, \cite{KS}]
\label{th-T-def}
Let $S\ni P$ be a germ of a two-dimensional quotient singularity.
The following are equivalent:
\begin{enumerate}
\item
$S\ni P$ is either Du Val or of type $\T$,
\item
$S\ni P$ has a $\QQ$-Gorenstein one-parameter smoothing,
\item
there is a terminal three-dimensional singularity $X\ni P$ and an
embedding $P\in S\subset X$ such that $S\subset X$ Cartier at $P$
and $(X,S)$ is plt.
\end{enumerate}
\end{theorem}

\section{Constructing semistable Mori conic bundles via deformations}
\label{S-def}
In this section all varieties are assumed to be analytic spaces.
\begin{definition}
A log (resp. $\T$) contraction $\var\colon S\to T\ni o$ with $\dim
T=1$ is called a \emph{log} (resp. $\T$) \emph{conic bundle}.
\end{definition}
\begin{theorem}
Let $\var \colon S \to T\ni o_T$ be a $\T$-conic bundle. There
exists a semistable Mori conic bundle $f\colon X\to Z\ni o_Z$ with
smooth base and embeddings
\[
\begin{diagram}[height=20pt,width=35pt]
S&\rInto &X
\\
\dTo^{\var}&&\dTo^f
\\
T&\rInto^{\upsilon} &Z
\end{diagram}
\]
such that $\upsilon(o_T)=o_Z$ and $(X,S)$ is plt.
\end{theorem}
We shall construct $X$ as an one-parameter family of
$T$-contractions.
\begin{proof}
We replace $S$ and $T$ with their compactifications so that $S$
and $T$ are projective, $T\simeq \PP^1$, and $\var \colon S \to T$
is smooth outside of $\var^{-1}(o_T)$. Let $P_i$ be singular
points of $S$.

Denote $\mathrm{Def}(S)$ (resp. $\mathrm{Def}(S,P_i)$) the base
space of the versal deformation of $S$ (resp. of the singularity
$S\ni P_i$). Let $s\colon \EuScript{S}\to \mathrm{Def}(S)$ be the
versal family. Thus we may assume that $S=s^{-1}(o)$ for some
$o\in \mathrm{Def}(S)$.

From \cite[Proposition 11.4]{KM1} we obtain the existence of the
diagram of morphisms of complex analytic spaces.
\[
\begin{diagram}[height=20pt,width=30pt]
\EuScript{S}&\rTo^{\scriptstyle\EuScript{F}}&\EuScript{T}=T\simeq\PP^1
\\
\dTo&&\dTo
\\
\mathrm{Def}(S)& \rTo& \mathrm{Def}(T)=\mt{pt}
\end{diagram}
\]
where $\EuScript{F}(S)=T$ and $\EuScript{F}|_S=\var$.

There is a natural (pull-back) morphism of germs of analytic
spaces
\begin{equation}
\label{eq-map-def}
\mathrm{Def}(S)\longrightarrow \prod_i \mathrm{Def}(S,P_i)
\end{equation}

The obstruction to globalize deformations in $\prod_i
\mathrm{Def}(S,P_i)$ lies in $R^2\var_*\Theta_S$, where
$\Theta_S=(\Omega^1_S)^\vee$, the tangent sheaf of $S$. Since
$R^2\var_*\Theta_S=0$, the map \eqref{eq-map-def} is smooth. In
particular, every deformation of singularities $S\ni P_i$ may be
globalized (cf. \cite[11.4.2]{KM1}).

By Theorem \xref{th-T-def} every singularity of $S$ admits a
$\QQ$-Gorenstein one-parameter smoothing. Therefore there is a
smoothing $g\colon X\to \varDelta\ni 0$, where $g^{-1}(0)=S$, $X$
is $\QQ$-Gorenstein and $\varDelta\subset \CC$ is a small disc. By
Inversion of Adjunction $(X,S)$ is plt and since $S$ is Cartier,
$X$ has at worst terminal singularities near $S$.

Put $Z=T\times \varDelta$ and let $f\colon X\to Z$ be the
projection. It is clear that $f|_S=\var\colon S\to T$. Therefore
$-K_X$ is $f$-ample near $S$. Shrinking $Z$ we get a Mori conic
bundle $f\colon X\to Z\ni o=(o_T,0)$.
\end{proof}

\section{Two-dimensional log contractions}
\label{sect-lcb}
\begin{notation}
\label{2-notation}
Let $\var\colon S\to T\ni o$ be a log contraction. We assume that
$S$ has at least one non-Du Val singularity. Let $\mu\colon \tilde
S\to S$ be a minimal resolution and let $\phi\colon \tilde S\to T$
be the composition map. Take an effective Cartier divisor $D$ on
$S$ such that $\Supp(D)=\var^{-1}(o)$ and $-D$ is $\var$-nef. For
example, in the case $\dim T=1$ we can put $D=\var^*(o)$ (the
scheme-theoretical fiber) while in the case $\dim T=2$ we can put
$D:=\var^*\var_*H-H$, where $H$ is a very ample divisor on $S$
such that $\var_*H$ is Cartier.

One can write the standard formula
\begin{equation}\label{eq-codiscr}
\mu^*K_S=K_{\tilde S}+\De,
\end{equation}
where $\De$ is an effective exceptional divisor, so-called,
\emph{codiscrepancy divisor}. Since the singularities of $S$ are
log terminal, $\up \De=0$. We also write $\mu^* D=\sum
l_iL_i+e_jE_j$, where the $E_i$ (resp. $L_i$) are
$\mu$-exceptional (resp. non-$\mu$-exceptional) components and
$l_i,e_i\in \NN$. Put $L=\sum l_iL_i$ and $E=e_jE_j$. Thus,
$D=\mu_*L$ and $\Supp(\De)\subset \sum E_i$.
\end{notation}

\begin{lemma}
\label{lemma-crit-amp}
Notation as above.
\begin{enumerate}
\item
$\Supp(L+E)$ is a tree of smooth rational curves;
\item
all the components of $L$ are $(-1)$-curves;
\item
$\De\cdot L_i<1$ for all $i$;
\item
if $\dim T=1$, then $L\cdot \De+2=\sum l_i$.
\end{enumerate}
\end{lemma}

\begin{proof}
{\rm (i)} is obvious because $\phi$ is flat in the case $\dim T=1$
and because $Z\ni o$ is a rational singularity in case $\dim T=2$.
By \eqref{eq-codiscr} we have
\[
0>\mu^*K_S\cdot L_i=K_{\tilde S}\cdot L_i+\De\cdot L_i=\De\cdot
L_i-2-L_i^2.
\]
Since $L_i^2<0$ and $\De\cdot L_i\ge 0$, we have $L_i^2=K_{\tilde
S}\cdot L_i=-1$ and $\De\cdot L_i<1$. This shows {\rm (ii)} and
{\rm (iii)}. For {\rm (iv)} we note that $\mu^*K_S\cdot L=-2$.
Thus,
\[
\sum l_i=-K_{\tilde S}\cdot L=-\mu^*K_S\cdot L+\De\cdot L=\De\cdot
L+2.
\]
\end{proof}

\begin{remark}
\label{rem-crit-amp}
\begin{enumerate}
\item
It is easy to see that the condition (iii) of
\xref{lemma-crit-amp} is also sufficient. Assume that conditions
\xref{2-notation} hold except for the ampleness of $-K_S$. If
$\De\cdot L_i<1$ for all $i$, then $\var$ is a log contraction,
i.e., $-K_S$ is ample.
\item
If $S$ has a unique non-Du Val point, then (iii) holds.
\end{enumerate}
\end{remark}

To describe log contractions we make use the weighted graph
language. By a \emph{weighted} graph $\Ga$ we mean a graph where
each vertex is given a weight $b_i\ge 1$. With a weighted graph
$\Ga=\{v_1,\dots,v_\varrho\}$ we associate a \emph{quadratic form}
by setting $v_i^2=-b_i$ and $v_i\cdot v_j$ for $i\neq j$ is equal
to the number of edges joining $v_i$ and $v_j$. We say that a
weighted graph $\Ga=\{v_1,\dots,v_\varrho\}$ is \emph{elliptic}
(resp. \emph{parabolic}) if the associated quadratic form has
signature $(0,\varrho)$ (resp. $(0,\varrho-1)$). Vertices with
$b_i=1$ will be referred to (and drawn) as \emph{black} vertices
those with $b_i\ge 2$ as \emph{white}. Weights $b_i=1$ and $b_i=2$
will be omitted.

By the \emph{blow-up of a vertex} $v_i$ we mean the following
transformation: the weight of the vertex $v_i$ increases by one,
that is, $b_i'=b_i+1$ and a new black vertex is added to the
graph, attached by an edge to the vertex $v_i'$. Similarly the
\emph{blow-up of an edge} $\{v_i,v_j\}$ is the following
transformation: the weight of the vertices $v_i$ and $v_j$
increase by one, the number of edges joining $v_i$ and $v_j$
decreases by one, and a new black vertex is added to the graph,
attached by edges to the vertices $v_i'$ and $v_j'$. The inverse
transformations are called \emph{contractions}. One can easily see
how the above transformations are related to blow-ups of curves on
a smooth surface.

We denote by $[a_1,\dots,a_r]$ a (weighted) chain and by $[p\mid
a_1,\dots,a_r \mid b_1,\dots,b_s \mid c_1,\dots,c_l]$ a fork $\Ga$
having the central vertex $v_0$ of weight $p$ so that
$\Ga\setminus \{v_0\}$ is a disjoined union of chains
$[a_1,\dots,a_r]$, $[b_1,\dots,b_s]$, and $[c_1,\dots,c_l]$, where
vertices corresponding $a_1$, $b_1$, and $c_1$ are adjacent to
$v_0$.

Now in notation \xref{2-notation} we denote by $\Ga(\var)$ the
dual graph of the fiber $\phi^{-1}(o)$. By (i) of Lemma
\xref{lemma-crit-amp}, $\Ga(\var)$ is a tree. Moreover, the graph
$\Ga(\var)$ is parabolic whenever $\dim T=1$ and elliptic whenever
$\dim T=2$.

\begin{lemma}
\label{l-graph-not}
The fork $[1 \mid a \mid b\mid c]$ is not elliptic for $a,b,c\ge
1$. The following graphs are parabolic \textup(and not
elliptic\textup{):}

{\rm (i)} chains $[1,1]$, $[1,2,\dots,2,1]$, $[2,1,2]$,

{\rm (ii)} the fork $[2 \mid 2 \mid 2\mid 2,\dots,2,1]$
\end{lemma}

\begin{corollary}
\label{cor-struct-log-c}
Let $D_i$ be irreducible components of $D$, then
\begin{enumerate}
\item
intersection points $D_i\cap D_j$ are singular and not Du Val,
\item
there are at most two singular points on every $D_i$,
\item
$(S,D_i)$ is plt near every Du Val point,
\item
$S$ has no Du Val singularities of type $D_n$ and $E_n$.
\end{enumerate}
\end{corollary}

\subsection*{Configuration of singular points}
\begin{lemma}
\label{lemma-not-plt-case}
Let $\var\colon S\to T\ni o$ be a $\T$-contraction and let $C$ be
a component of $D$. Assume that $C$ contains exactly two singular
points of $S$ and they are of type $\T$ \textup(not Du
Val\textup). Then $(S,C)$ is plt.
\end{lemma}
\begin{proof}
If $C\neq \Supp(D)$, then $C$ is contractible over $T$, i.e.,
there is a decomposition $S\to T'\to T$ such that $\var'\colon
S\to T'$ is birational and $C$ is the only $\var'$-exceptional
divisor. Replacing $T$ with $T'$ we may assume that $C= \Supp(D)$.

Let $P_1$, $P_2$ be singular points. Assume that $(S,C)$ is not
plt near $P_1$. We claim that $(S,C)$ is plt near $P_2$. Indeed,
take two general hyperplane sections $H_1$ and $H_2$ passing
through $P_1$ and $P_2$. For some $0<\ep'\ll\ep \ll 1$ the divisor
$-(K_S+(1-\ep')C+\ep H_1+\ep H_2)$ is $\var$-ample and the pair
$(S,(1-\ep')C+\ep H_1+\ep H_2)$ is not lc at $P_1$ and $P_2$. This
contradicts Connectedness Lemma \cite[5.7]{Sh}.

Thus $(S,C)$ is plt near $P_2$ and $\Ga(\var)$ has the form
\begin{equation}
\label{eq-non-plt-graph}
\begin{diagram}[abut,size=12pt]
\stackrel{b_1}{\circ}&\lLine&\cdots&\stackrel{b_k}{\circ} &\cdots&
\lLine&\stackrel{b_\varrho}{\circ}
\\
&&&\dLine
\\
&&&\bullet&\lLine&\stackunder{c_1}{\circ}&\lLine&
\cdots&\lLine&\stackunder{c_l}{\circ}&
\end{diagram}
\end{equation}
where $[b_1,\dots, b_\varrho]$ and $[c_1,\dots,c_l]$ are
$\T$-chains (i.e. $\Gamma(\var)=[b_k\mid b_{k-1},\dots,b_1\mid
b_{k+1},\dots,b_\varrho\mid c_1,\dots,c_l]$). Let $\alpha_1'$ and
$\alpha_k$ be log discrepancies of the vertices corresponding to
$c_1$ and $b_k$, respectively. By Lemma \xref{lemma-crit-amp} we
have $\alpha_1'+\alpha_k>1$. Let
\[
\Ga(\var)=\Ga_0\to \Ga_1\to\cdots\to \Ga_r=\Ga'
\]
be the sequence of contractions of black vertices. If $b_k=2$,
then $\Ga_1$ contains the fork $[1\mid b_{k-1}\mid b_{k+1}\mid
c_1-1]$. This contradicts Lemma \xref{l-graph-not}. Therefore
$b_k\ge 3$. The same arguments show that in each graph $\Ga_i$ the
central vertex (corresponding to $b_k$) is not black. Since
$[b_1,\dots,b_\varrho]$ is a chain of type $\T$, we may assume
that $b_\varrho=2$ and $b_1\ge 3$. Thus
\[
\Ga'=[b_k-s\mid b_{k-1},\dots,b_1\mid b_{k+1},\dots,b_\varrho\mid
c_s-1,c_{s+1},\dots,c_l],
\]
where $s\ge 1$, $b_k-s\ge 2$, and $c_s-1\ge 2$. Clearly, $\Ga'$ is
elliptic and log terminal (because $-(K_{\tilde S}+\De)$ is nef
and big over $T$). Now we use the classification of
two-dimensional log terminal singularities (see, e.g., \cite[Ch.
3]{Ut}). Since $b_1>2$, we have $[c_s-1,\dots,c_l]=n/q$, where
$1\le q<n$, $\gcd(n,q)=1$ and $n=2,3,4$, or $5$. So,
$[c_s-1,\dots,c_l]=[n]$, $[{2,\dots,2}]$, $[2,3]$, or $[3,2]$.

Assume that $s=1$. Then $[c_1,\dots,c_l]=[4]$ or $[3,3]$ (see
Proposition \xref{prop-T-ind}), $k=\varrho-1$, and
$[b_1,\dots,b_{k-1}]=n'/q'$, where $1\le q'<n'-1$, $\gcd(n',q')=1$
and $n'=3,4$, or $5$. Further, $\alpha_1'=1/2$. Easy computations
(see \cite[(3.1.3)]{Ut}) show that $2/b_k> \alpha_k>1/2$. Thus
$b_k= 3$. We get only one possibility
$[b_1,\dots,b_k,\dots,b_\varrho]=[4,3,2]$. But then
$\alpha_k=1/5$, a contradiction.

Assume that $s>1$. Then $c_1=\cdots=c_{s-1}=2$. Hence, $c_l\ge 3$
and $[c_1,\dots,c_l]= [2,\dots,2,n+1]$, or $[2,\dots,2,3,3]$.
Since $[c_1,\dots,c_l]$ is a $\T$-chain, $n\ge 4$ and
$[c_1,\dots,c_l]$= $[2,5]$, or $[2,2,6]$. As above we get
$\alpha_1'\le 1/3$, $2/b_k>\alpha_k>2/3$, and $b_k=2$, a
contradiction.
\end{proof}

\begin{lemma}
Notation as above. Let $D_1$ and $D_2$ be two components of $D$
such that
\begin{enumerate}
\item
$D_1\cap D_2\neq\emptyset$,
\item
there are $\T$-points $P_i\in D_i$, $P_i\notin D_1\cap D_2$.
\end{enumerate}
Then $K_S+D_1+D_2$ is lc.
\end{lemma}
\begin{proof}
Assume the converse. By the previous lemma, $\Ga(\var)$ contains a
subgraph $\Ga$ of the form
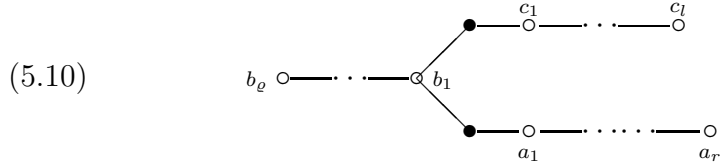
\begin{equation}
\label{eq-plt-2-graph}
\begin{diagram}[abut,size=10pt]
&&&& &&&\bullet&\lLine&\stackrel{c_1}\circ&\lLine&\cdots&
\lLine&\stackrel{c_l}\circ
\\
&&&& &&\ruLine&&&&&
\\
\scriptstyle{b_\varrho}&\circ&\lLine&\cdots&
\lLine&\circ&\scriptstyle{b_1}&&&& &
\\
&&&& &&\rdLine&&&&&
\\
&&&& &&&\bullet&\lLine&\stackunder{a_1}\circ&\lLine&\cdots&\cdots&
\lLine&\stackunder{a_r}\circ
\end{diagram}
\end{equation}
where $\varrho\ge 2$ and $[a_1,\dots,a_r]$, $[b_1,\dots,
b_\varrho]$, $[c_1,\dots,c_l]$ are $\T$-chains.

Note that $b_1\ge 3$. By Corollary \xref{cor-alpha-alpha} and
Lemma \xref{lemma-crit-amp} we have $a_1=c_1=2$. Take $s,m\ge 1$
so that
\[
a_1=\cdots=a_{s}=2,\ a_{s+1}>2,\ c_1=\cdots=c_{m}=2,\ c_{m+1}>2.
\]
Contracting black vertices successively we get the following log
terminal graph
\[
\Gamma'=[b_1-s-m-2\mid b_2,\dots,b_\varrho\mid
a_{s+1}-1,\dots,a_r\mid c_{m+1}-1,\dots,c_l].
\]
By Proposition \xref{prop-T-ind} we have $\sum a_i= 2-\beta+3r$,
where $\beta$ is the number of vertices of the corresponding chain
\eqref{eq-graph-basic} with $\iota=2$ ($\beta$ coincides with
$\beta(n/q)$ introduced in \S \xref{sect-T} but we do not need
this fact). Since $\beta+s\le r$,
\[
a_{s+1}+\cdots+a_r=2-\beta+3r-2s\ge 2+r+\beta\ge 5.
\]
Similarly, $c_{m+1}+\cdots+c_l\ge 5$. Therefore,
$\varrho=b_\varrho=2$ and we may assume that
$[a_{s+1}-1,\dots,a_r]=3/q$. On the other hand, $a_r\ge 3$, so
$[a_{s+1}-1,\dots,a_r]=[3]$ and $[a_1,\dots,a_r]=[2,\dots,2,4]$, a
contradiction.
\end{proof}

\begin{corollary}
Let $\var\colon S\to T\ni o$ be a $\T$-contraction.
\begin{enumerate}
\item
If $S$ has exactly one non-Du Val point $P$, then all the
components of $D$ pass through $P$.
\item
If $S$ has more than one non-Du Val points, then $(S,D_i)$ is plt
for any component $D_i$ of $D$ containing two non-Du Val points.
\end{enumerate}
\end{corollary}

Now Theorem \xref{ge} is a consequence of Propositions
\xref{prop-compl-ext} and \xref{prop-1-compl} below.

\begin{proposition}
\label{prop-1-compl}
Let $\var\colon S\to T\ni o$ be a $\T$-contraction. Then $K_S$ is
$1$-complementary.
\end{proposition}

\begin{pusto}
\label{notation-3}
To begin with, assume that $\var\colon S\to T\ni o$ is an
\emph{arbitrary} log contraction. We apply the technique developed
in \cite{PSh}. Take $\delta$ so that $(S,\delta D)$ is maximally
log canonical.
\end{pusto}

Note that on $S$ the LMMP works with respect to any divisor $G$.
Indeed, there is a boundary $F$ such that $K_X+F$ is klt,
numerically trivial, and the components of $F$ generate
$\Pic(S)\otimes \QQ$. Then $G$-MMP is equivalent to $K_S+F+\ep
G'$-MMP for $0<\ep\ll 1$ and suitable $G'\qq G$.

\begin{lemma}
\label{lemma-plt-case}
Assume that $(S,\delta D)$ is plt. Then $K_S$ is
$1$-complementary.
\end{lemma}
\begin{proof}
Since $(S,\delta D)$ is maximally log canonical, $\down{\delta
D}\neq 0$. Put $C=\down{\delta D}$ and $B=\delta D-C$. By
Connectedness Lemma \cite[5.7]{Sh}, $C$ is an irreducible curve.
By Corollary \xref{cor-struct-log-c}, $\Diff_C(\delta B)$ is
supported in two points, say $P_1$ and $P_2$. Run $-(K_S+C)$-MMP
over $T$:\quad $\psi\colon S\to \bar S$. Since $-K_S$ is
$\psi$-ample, $C$ is not contracted. Since $-(K_S+\delta D)$ is
$\psi$-ample, the plt property of $(S,\delta D)$ is preserved. We
get a plt model $(\bar X,\bar C)$ such that $-(K_{\bar S}+\bar C)$
is nef over $T$. By the above, $\Diff_{\bar C}(\delta \bar B)$ is
supported in two points. Hence $K_{\bar S}+\Diff_{\bar C}$ is
$1$-complementary (see \cite[5.2]{Sh}). Since $-(K_{\bar S}+\bar
C)$ is nef and big over $T$, this complement can be extended to
$\bar S$ (see \cite[Prop. 4.4.1]{Lect}). By \cite[4.3]{Lect} this
gives us an $1$-complement of $K_S+C$.
\end{proof}

\begin{pusto}
\label{notation-4}
If $(S,\delta D)$ is not plt, there is an \emph{inductive blowup}
of $(S,\delta D)$. By definition it is a birational extraction
$\sigma\colon S'\to S$ with irreducible exceptional divisor $E$
satisfying the following properties:
\begin{enumerate}
\item
$K_{S'}+E+\delta D'=\sigma^*(K_S+\delta D)$ is log canonical,
where $D'$ is the proper transform of $D$,
\item
$(S', E)$ is plt.
\end{enumerate}
Since the minimal resolution $\mu \colon \tilde S\to S$ is a log
resolution of $(S,D)$, we may assume that $\mu$ factors through
$\sigma$ (see \cite[Proof of 3.1.4]{Lect}). Then $K_{S'}+\alpha
E=\sigma^*K_S$, where $\alpha\ge 0$ and $-(K_{S'}+\alpha E)$ is
nef over $T$. As in the proof of Lemma \xref{lemma-plt-case} we
run $-(K_{S'}+E)$-MMP over $T$. In this case again $E$ cannot be
contracted and we get a model $(\bar S,\bar E)$ such that
$-(K_{\bar S}+\bar E)$ is nef over $T$:
\[
\begin{diagram}[height=20pt]
&&S'&&
\\
&\ldTo^{\sigma}&&\rdTo^{\psi}&
\\
S&&&&\bar S
\\
&\rdTo^{\var}&&\ldTo^{\bar \var}&
\\
&&T&&
\end{diagram}
\]
Denote by $N$ the exceptional divisor of $\psi$ and let
$V:=\Sing(S')\cap E$. If $-(K_{S'}+E)$ is nef over $T$, we put
$\psi=\mt{id}$ and $N=\emptyset$. Clearly all singular points
$\bar P_1,\dots,\bar P_r$ of $\bar S$ lying on $\bar E$ are
contained in $\psi(V)\cup \psi(N)$. If $r\le 2$, then $K_{\bar
E}+\Diff_{\bar E}$ is $1$-complementary (see \cite{Sh}). Since
$-(K_{\bar S}+\bar E)$ is nef over $T$, this complement can be
extended to $\bar S$ (see \cite[Prop. 4.4.1]{Lect}). By
\cite[4.3]{Lect} this gives us an $1$-complement of $K_S$.
\end{pusto}

\begin{lemma}
Assume that $S$ has a unique non-Du Val point and this point which
is a cyclic quotient. Then $K_S$ is $1$-complementary.
\end{lemma}
\begin{proof}
We may assume that $(S,\delta D)$ is not plt. Since
$P:=\sigma(E)\in S$ is a cyclic quotient singularity, $V$ contains
at most two points. Indeed, $-K_{S'}$ is $\psi$-ample and $S'$ has
at worst Du Val singularities outside of $\Sing(\bar S)\cap \bar
E$. By \cite[3.38]{KM2} discrepancies of all divisors of $\bar S$
over $\psi(E)$ are strictly positive. Hence $\bar S$ is smooth at
points of $\psi(N)$, $\bar P_1,\dots,\bar P_r\in \psi(V)$ and
$r\le 2$. By the above $K_S$ is $1$-complementary.
\end{proof}

\begin{proof}[Proof of Proposition {\xref{prop-1-compl}}]
Again $V$ contains at most two points. If $K_S$ is not
$1$-complementary, then $r\ge 3$. Take $\bar P\in \psi(N)\setminus
\psi(V)$ and let $N_{0}=\psi^{-1}(\bar P)$. Then the point
$N_{0}\cap E\in S$ is smooth and $N_{0}$ contains at least one
non-Du Val point of $S'$. If $V$ is two points, then by Corollary
\xref{cor-struct-log-c} we get a subgraph
\eqref{eq-non-plt-graph}, a contradiction.

Assume that $V$ is one point. Then there are two points $\bar P,
\bar P_1\in \psi(N)\setminus \psi(V)$ and similarly by Lemma
\xref{lemma-not-plt-case} we get a subgraph
\eqref{eq-plt-2-graph}, a contradiction.

Finally, assume that $S'$ is smooth along $E$. Then $E$ is a
$(-4)$-curve. Hence $\Ga(\var)$ contains a fork $[4\mid 1,b_1\mid
1,b_2\mid 1,b_3]$. Such a graph cannot be elliptic.
\end{proof}

\subsection*{Examples}
\input border.tex

\begin{proposition}
\label{prop-indx-2}
Let $\var\colon S\to T$ be a two-dimensional log conic bundle of
index two. Then $\Ga(\var)$ is one if the following:
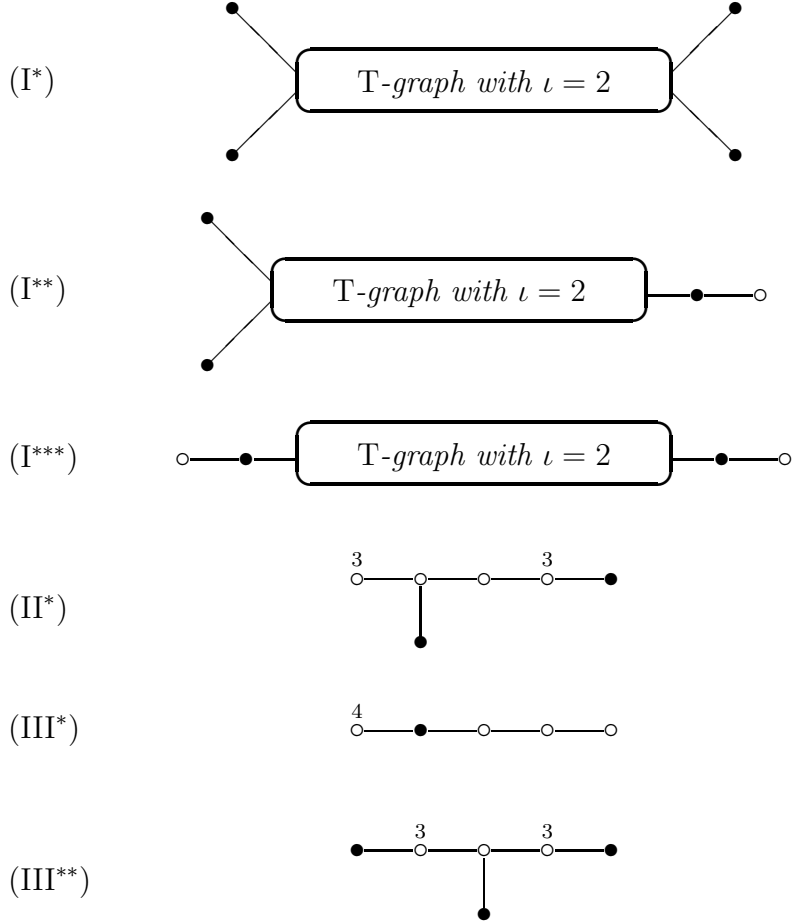
\begin{equation*}
\begin{diagram}[abut,size=12pt]
\bullet&&&&\bullet
\\
&\luLine&&\ruLine&
\\
&& \hspace{7pt}
\begin{tabular}{c}
\setborder[5cm,0.2cm,0.2cm]\border{\center{$\T$-graph with
$\iota=2$}}
\end{tabular}
\hspace{7pt}&&
\\
&\ldLine&&\rdLine&
\\
\bullet&&&&\bullet
\end{diagram}
\leqno (\mathrm I^*)
\end{equation*}

\begin{equation*}
\begin{diagram}[abut,size=12pt]
\bullet&&&&&&
\\
&\luLine&&&&&
\\
&& \hspace{7pt}
\begin{tabular}{c}
\setborder[5cm,0.2cm,0.2cm]\border{\center{$\T$-graph with
$\iota=2$}}
\end{tabular}
\hspace{-6pt} &\lLine&\bullet&\lLine&\circ
\\
&\ldLine&&&&&
\\
\bullet&&&&&&
\end{diagram}
\leqno (\mathrm I^{**})
\end{equation*}

\begin{equation*}
\begin{diagram}[abut,size=12pt]
\circ&\lLine&\bullet&\lLine& \hspace{-6pt}
\begin{tabular}{c}
\setborder[5cm,0.2cm,0.2cm]\border{\center{$\T$-graph with
$\iota=2$}}
\end{tabular}
\hspace{-6pt} & \lLine&\bullet&\lLine&\circ
\end{diagram}
\leqno (\mathrm I^{***})
\end{equation*}

\begin{equation*}
\begin{diagram}[abut,size=12pt]
\stackrel{3}{\circ}&\lLine&\circ&\lLine&
\circ&\lLine&\stackrel{3}{\circ}&\lLine&\bullet
\\
&&\dLine&&&&&&
\\
&&\bullet&&&&&&
\end{diagram}
\leqno (\mathrm{II}^*)
\end{equation*}

\begin{equation*}
\begin{diagram}[abut,size=12pt]
\stackrel{4}{\circ}&\lLine&\bullet&\lLine&\circ&\lLine&\circ&\lLine&\circ
\end{diagram}
\leqno (\mathrm{III}^*)
\end{equation*}

\begin{equation*}
\begin{diagram}[abut,size=12pt]
\bullet&\lLine&\stackrel{3}{\circ}&\lLine&\circ&
\lLine&\stackrel{3}{\circ}&\lLine&\bullet
\\
&&&&\dLine&&&&
\\
&&&&\bullet&&&&
\end{diagram}
\leqno (\mathrm{III}^{**})
\end{equation*}
\end{proposition}

Our notation can be explained as follows. Consider a general
member $B\in |-K_S|$ and let $S'\to S$ be a double covering
branched along $B$. Then $K_{S'}=0$, $S'\to T$ is an elliptic
fibration, and $S'$ has only Du Val singularities (cf. \cite[\S
3]{P-sb}). If $\tilde S'$ is the minimal resolution, then the
central fiber of the composition map $\tilde S'\to T$ is Kodaira's
singular fiber.

\begin{proof}
For any index two log terminal point all log discrepancies of the
minimal resolution are equal to $1/2$. By Lemma
\xref{lemma-crit-amp} there is at most one non-Du Val point on
each component of $D$. By Corollary \xref{cor-struct-log-c} there
is exactly one non-Du Val point $P$ on $S$ and all the components
of $D$ pass through $P$. Again using Lemma \xref{lemma-crit-amp}
we have $\sum l_i=4$, so the graph $\Ga(\var)$ has at most four
black vertices. Now the classification follows by Lemma
\xref{l-graph-not}.
\end{proof}

\begin{example}
\label{ex-2-3}
The following graphs gives us examples of $\T$-conic bundles with
two and three non-Du Val points.
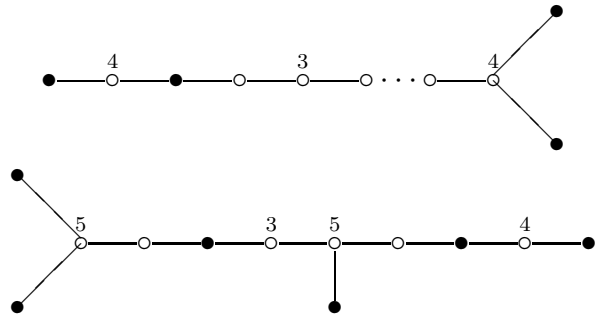
\[
\begin{diagram}[abut,size=12pt]
&&&&&& &&&&&&&&&&\bullet
\\
&&&&&& &&&&&&&&&\ruLine&
\\
\bullet&\lLine&\stackrel{4}{\circ}&\lLine&\bullet&\lLine&\circ&\lLine
&\stackrel{3}{\circ}&\lLine&\stackrel{}{\circ}&
\cdots&\stackrel{}{\circ}&\lLine&\stackrel{4}{\circ}&
\\
&&&&&& &&&&&&&&&\rdLine&
\\
&&&&&& &&&&&&&&&&\bullet
\end{diagram}
\]
\begin{equation*}
\begin{diagram}[abut,size=12pt]
\bullet&&& &&& &&& &&& &&& &&& &
\\
&\luLine&& &&& &&& &&& &&& &&& &
\\
&&\stackrel5\circ&\lLine &\circ&\lLine&\bullet
&\lLine&\stackrel3\circ&\lLine &\stackrel5\circ&\lLine&\circ
&\lLine&\bullet&\lLine &\stackrel4\circ&\lLine&\bullet &
\\
&\ldLine&& &&& &&& &\dLine&& &&& &&& &
\\
\bullet&&& &&& &&& &\bullet&& &&& &&& &
\end{diagram}
\end{equation*}
\end{example}

\begin{proposition}
For any $\T$-singularity $\frac1n(1,q)$ there is a $\T$-conic
bundle having exactly one singular point which is of type
$\frac1n(1,q)$.
\end{proposition}
\begin{proof}
One can start with graph $(\mathrm{I}^*)$ and run the construction
below. By Proposition \xref{prop-T-ind} on each step we have only
singularities of type $\T$.
\end{proof}

\begin{construction}
\label{constr-1}
Let $\var$ be a log conic bundle with a unique singular point that
is of type $[b_1,\dots,b_\varrho]$. Assume that in $\Ga(\var)$
there is a black vertex adjacent to the end $b_1$, i.e.,
$\Ga(\var)$ contains a string $[1,b_1,\dots,b_\varrho]$, where
$[b_1,\dots,b_\varrho]$ corresponds to singular point. Blowing-up
the ends $1$ and $b_\varrho$, we get a graph $\Gamma'$ containing
a string $[1,2,b_1,\dots,b_\varrho+1,2]$. By Remark
\xref{rem-crit-amp}, $\Gamma'$ corresponds to a log conic bundle
(i.e., the anticanonical divisor is ample).
\end{construction}

\begin{remark}
\begin{enumerate}
\item
Construction \xref{constr-1} provides only one series of
$\T$-conic bundles with a unique singular point. For example, the
following $\T$-conic bundle cannot be obtained by this way.
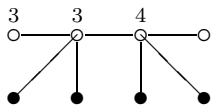
\begin{equation*}
\begin{diagram}[abut,size=12pt]
\stackrel3\circ& \lLine& \stackrel3\circ& \lLine& \stackrel4\circ&
\lLine& \circ
\\
&\ldLine&\dLine&&\dLine&\rdLine&
\\
\bullet&&\bullet&&\bullet&&\bullet
\end{diagram}
\end{equation*}
\item
Similar to \xref{constr-1} one can obtain infinite series of
$\T$-conic bundles with two and three singular points starting
with Example \xref{ex-2-3}.
\end{enumerate}
\end{remark}

\section{The case of irreducible central fiber}
The following lemma shows that case of relative Picard number one
is most important.
\begin{lemma}
Let $f\colon X\to Z\ni o$ be a Fano-Mori contraction such that the
dimension of fibers is at most one. There is a Fano-Mori
contraction $f'\colon X'\to Z$ with the same property and a
birational map $h\colon X\dashrightarrow X'$ over $Z$ such that
$h$ is a morphism outside of $f^{-1}(o)$, $X'$ is $\QQ$-factorial
and $\rho(X/Z)=1$. In particular, $f'^{-1}(L)$ is irreducible for
any irreducible divisor $L\subset Z$. If furthermore $(X,f^*T)$ is
dlt for some effective Cartier divisor $T$, then we can take $X'$
so that $(X',f'^*(T))$ is dlt.
\end{lemma}
\begin{proof}
Let $q\colon X^{\mathrm q}\to X$ be a $\QQ$-factorial
modification. Thus $X^{\mathrm q}$ has only terminal
$\QQ$-factorial singularities, $K_{X^{\mathrm q}}=q^*K_X$, and $q$
is a small birational contraction. Run MMP over $Z$. We get $X'$
as above.

To show the second statement we construct $\QQ$-factorialization
$q\colon X^{\mathrm q}\to X$ of $(X,S)$. Then $(X^{\mathrm
q},S^{\mathrm q})$ is dlt, where $S^{\mathrm q}=q^*S$. Then we
just note that $K_{X^{\mathrm q}}$-MMP is the same as
$K_{X^{\mathrm q}}+S^{\mathrm q}$-MMP.
\end{proof}

In analytic situation $\rho(X/Z)=1$ implies that the fiber
$f^{-1}(o)$ is irreducible. Now we classify semistable Mori conic
bundles with irreducible central fiber.

\begin{proposition}[cf. {\cite[Th. 2.5]{P-rb}}]
\label{prop-irred-fib}
Let $\var \colon S\to T\ni o$ be a $\T$-conic bundle having at
least one non-Du Val point. Assume that the fiber $D$ is
irreducible. Then $\var$ is of type $(\mathrm{III}^*)$ of
\xref{prop-indx-2}.
\end{proposition}

\begin{corollary}
Let $f\colon X\to Z\ni o$ be a semistable Mori conic bundle such
that $f^{-1}(o)$ is irreducible. Then $X$ has exactly one
non-Gorenstein point which is of index $2$, see \cite[\S 3]{P-sb}.
\end{corollary}

\begin{proof}
If $K_S+C$ is plt, then $S$ has two singularities of types
$\frac1n(1,q)$ and $\frac1n(1, n-q)$ (see \cite[Th. 2.5]{P-rb}).
If they are of type $\T$, then
\[
(q+1)^2\equiv 0\mod n, \qquad (n-q+1)^2\equiv 0\mod n.
\]
This gives us $4\equiv 0\mod n$. Since the singularities of $S$
are worse than Du Val, $n=4$. We get case $(\mathrm{III}^*)$.

Now we consider the case when $K_S+C$ is not plt. By Lemma
\xref{lemma-not-plt-case} and Corollary \xref{cor-struct-log-c}
the surface $S$ has exactly one non-Du Val point and at most one
Du Val point $Q$. Moreover, $S\ni Q$ is of type $A_n$ and $K_S+C$
is plt near $Q$. Thus $\Ga(\var)$ has the following form
\begin{equation}
\label{eq-graph-o}
\begin{diagram}[abut,size=12pt]
\stackrel{b_1}{\circ}&\lLine&\stackrel{b_2}{\circ}&
\lLine&\cdots&\lLine&\stackrel{b_r}{\circ}&
\lLine&\cdots&\lLine&\stackrel{b_\varrho}{\circ}
\\
&&&&&&\dLine&&&&
\\
&&&&&&\bullet&&&&
\\
&&&&&&\dLine&&&&
\\
{\circ}&\lLine& {\circ}&\lLine&\cdots&\lLine& {\circ}&&&&
\end{diagram}
\end{equation}
where $r\neq 1$, $r\neq \varrho$. Contracting black vertices
successively, on some step we get a subgraph
\begin{equation}
\label{eq-graph-par-line}
\begin{diagram}[abut,size=12pt]
\stackrel{b_1}{\circ}&\lLine&\cdots&\lLine&
\stackrel{b_{r-1}}{\circ}&\lLine&\bullet&\lLine&
\stackrel{b_{r+1}}{\circ}&\lLine&\cdots&\lLine&
\stackrel{b_\varrho}{\circ}
\end{diagram}
\end{equation}

\begin{lemma}
If the graph \eqref{eq-graph-par-line} is parabolic, then
\begin{equation}
\label{eq-cond}
\sum_{i=1}^{r-1}(b_i-1)=\sum_{j=r+1}^{\varrho} (b_j-1)= \varrho-2.
\end{equation}
In particular, $r\neq 1, \varrho$.
\end{lemma}
Apply the procedure described in Proposition \xref{prop-T-ind} to
$[b_1,\dots,b_r,\dots,b_\varrho]$. Each step preserves relation
\eqref{eq-cond}. At the end we get a chain
$[b_1',\dots,b_{r'}',\dots,b_{\varrho'}']$ as in
\eqref{eq-graph-basic}. Relation \eqref{eq-cond} give us
$r'=\varrho'-r'+1=\varrho'-2$, i.e., $r'=3$ and $\varrho'=5$.
Hence, in \eqref{eq-graph-par-line} we have $b_{r-1}=b_{r+2}=2$.
This contradicts Lemma \xref{l-graph-not}.
\end{proof}

\section{The existence of semistable $3$-fold flips}
\label{sect-flips}
\begin{theorem}
\label{th-flips}
Let $f\colon X\to Z$ be a semistable three-dimensional flipping
contraction. Assume that $f$ is extremal \textup(i.e., $X$ is
$\QQ$-factorial and $\rho(X/Z)=1$\textup). Then the flip of $f$
exists.
\end{theorem}
\begin{proof}[Sketch of the proof
\textup(see {\cite[8.5, 8.7]{Kaw}}\textup)] The existence of the
flip is equivalent to the finite generation of the
$\OOO_Z$-algebra $\RRR_Z(K_Z):=\oplus_{m\ge 0} \OOO_Z(mK_Z)$, see
\cite[Lemma 3.1]{Kaw}. By Theorem \xref{ge} there is $L=2F\in
|-2K_X|$ such that $K_X+f^*T+\frac12 L$ is lc. Since $f$ is an
extremal contraction and $K_X+f^*T+\frac12 L$ is numerically
trivial, one can see that $K_Z+T+\frac12L_Z$ is also lc, where
$L_Z:=f_*L\in |-2K_Z|$. Therefore, the same holds for a general
member $L_Z\in |-2K_Z|$ which is reduced and irreducible. As in
\cite[\S 8]{Kaw}, consider a double covering $\pi\colon Z'\to Z$
ramified along $L_Z$. Then $K_{Z'}+\pi^*T= \pi^*(K_Z+T+\frac12
L_Z)$ is lc and Cartier. Since $\pi^*T$ also is a Cartier divisor,
$Z'$ has at worst a canonical Gorenstein singularity at
$o':=\pi^{-1}(o)$. Put $L_{Z'}:=\pi^{*}(L_Z)_{\red}$. By
\cite[3.2]{Kaw} the finite generation of algebras $\RRR_{Z}(K_Z)$
and $\RRR_{Z'}(K_Z'-L_Z')$ is equivalent. Finally, the last
algebra is finitely generated by \cite[6.1]{Kaw} (see also
\cite{Kollar}, \cite[\S 6]{Ut}, \cite[\S 6]{KM2}).
\end{proof}

Note that in our case the finite generation of
$\RRR_{Z'}(K_Z'-L_Z')$ is much easier because the presence of a
Cartier divisor $\pi^*T$ such that $K_{Z}'+\pi^*T$ is lc.

\begin{corollary}[\cite{KM1}]
Let $f\colon X\to Z\ni o$ be a semistable birational contraction
with fibers of dimension at most one \textup(either flipping or
divisorial of type $2$-$1$\textup) and let $T$ be a general
hyperplane section through $o$. Then $T\ni o$ is a cyclic quotient
singularity. If furthermore $f$ is divisorial, then $T\ni o$ is of
type $\T$.
\end{corollary}

\begin{proof}
By Theorem \xref{ge} the pair $(Z,T+f(F))$ is log canonical. By
Adjunction so is $(T,\Diff_T(f(F)))$. Moreover,
$K_T+\Diff_T(f(F))\sim 0$. By the classification of log terminal
singularities with a reduced boundary $T\ni o$ is a cyclic
quotient singularity.
\end{proof}

\end{document}

%% file: border.tex
\ifx\BORDER\relax \else\let\BORDER=\relax\fi 
%
%
%
%
%
%

%


\edef\borderelm{\the\catcode`\_} 
\catcode`\_=11



\def\hborder#1#2#3#4
   {\hbox to#1
	{\setbox0=\hbox{\borderelm#3}%
	 \borderelm#2%
	 \hskip0pt minus1fil %
	 \copy0\kern-.5\wd0 %
	 \cleaders\copy0\hskip0pt plus1fil %
	 \kern-.5\wd0\copy0 %
	 \hskip0pt minus1fil %
	 \borderelm#4%
	}%
   }


\def\border#1#2
   {\vbox to#1
        {\offinterlineskip\boxmaxdepth=\maxdimen
	 \dimen2=#2
	 \setbox0=\hbox to\dimen2{\borderelm3\hss\borderelm4}
	 \dimen0=\ht0 \advance\dimen0 by\dp0
	 \hborder{\dimen2}012
	 \vskip0pt minus1fil
	 \copy0\kern-.5\dimen0
	 \cleaders\copy0\vskip0pt plus1fil
	 \kern-.5\dimen0\copy0
	 \vskip0pt minus1fil
	 \hborder{\dimen2}567
	 \kern0sp
	}%
   }


\def\hxborder#1#2#3#4
   {\hbox to#1
	{\setbox0=\hbox{\borderelm#3}%
	 \borderelm#2%
	 \xleaders\copy0\hskip0pt plus1fil minus1fil %
	 \borderelm#4%
	}%
   }


\def\xborder#1#2
   {\vbox to#1
        {\offinterlineskip\boxmaxdepth=\maxdimen
	 \dimen2=#2
	 \setbox0=\hbox to\dimen2{\borderelm3\hss\borderelm4}
	 \hxborder{\diemn2}012
	 \xleaders\copy0\vskip0pt plus1fil minus1fil
	 \hxborder{\dimen2}567
	 \kern0sp
	}%
   }



\def\hminusborder#1#2#3#4
   {\hbox
	{\dimen0=#1\count2=\dimen0 %
	 \setbox0=\hbox{\borderelm#2}\advance\count2 by-\wd0 %
	 \setbox0=\hbox{\borderelm#4}\advance\count2 by-\wd0 %
	 \setbox0=\hbox{\borderelm#3}%
	 \advance\count2 by100 %
	 \divide\count2 by\wd0 %
	 \borderelm#2%
	 \cleaders\copy0\hskip\count2\wd0 %
	 \borderelm#4%
	}%
   }

\def\minus_border#1#2#3
   {\vbox
        {\offinterlineskip\boxmaxdepth=\maxdimen
	 \dimen2=#2\relax
	 \dimen0=#1\count2=\dimen0
	 \setbox0=\hminusborder{\dimen2}012
	 \advance\count2 by-\ht0 \advance\count2 by-\dp0
	 \setbox2=\hminusborder{\dimen2}567
	 \advance\count2 by-\ht2 \advance\count2 by-\dp2
	 \if.#3.\setbox4=\hbox to\wd0{\borderelm3\hss\borderelm4}
	 \else  \setbox4=\hminusborder{\dimen2}384
	 \fi
	 \advance\count2 by100
	 \dimen0=\ht4\advance\dimen0 by\dp4
	 \divide\count2 by\dimen0
	 \box0
	 \cleaders\copy4\vskip\count2\dimen0
	 \box2
	 \kern0sp
	}%
   }

\def\minusborder#1#2{\minus_border{#1}{#2}{}}

\def\minusgray#1#2{\minus_border{#1}{#2}{1}}


\def\hplusborder#1#2#3#4
   {\hbox
	{\dimen0=#1\count2=\dimen0 %
	 \setbox0=\hbox{\borderelm#2}\advance\count2 by-\wd0 %
	 \setbox0=\hbox{\borderelm#4}\advance\count2 by-\wd0 %
	 \setbox0=\hbox{\borderelm#3}%
	 \advance\count2 by\wd0\advance\count2 by-100 %
	 \divide\count2 by\wd0 %
	 \borderelm#2%
	 \cleaders\copy0\hskip\count2\wd0 %
	 \borderelm#4%
	}%
   }

\def\plus_border#1#2#3
   {\vbox
        {\offinterlineskip\boxmaxdepth=\maxdimen
	 \dimen2=#2\relax
	 \dimen0=#1\relax\count2=\dimen0
	 \setbox0=\hplusborder{\dimen2}012
	 \advance\count2 by-\ht0 \advance\count2 by-\dp0
	 \setbox2=\hplusborder{\dimen2}567
	 \advance\count2 by-\ht2 \advance\count2 by-\dp2
	 \if.#3.\setbox4=\hbox to\wd0{\borderelm3\hss\borderelm4}
	 \else\setbox4=\hplusborder{\dimen2}384
	 \fi
	 \dimen0=\ht4\advance\dimen0 by\dp4
	 \advance\count2 by\dimen0\advance\count2 by-100
	 \divide\count2 by\dimen0
	 \box0
	 \cleaders\copy4\vskip\count2\dimen0
	 \box2
	 \kern0sp
	}%
   }

\def\plusborder#1#2{\plus_border{#1}{#2}{}}

\def\plusgray#1#2{\plus_border{#1}{#2}{1}}


\ifx\manfnt\UNDEFINED
\font\manfnt=manfnt
\fi


\def\roundcorners#1
   {{\manfnt
     \setbox2=\hbox{\char'44}
     \ifcase#1
	\setbox0=\hbox{d}
	\hbox to\ht0{\vbox to\fontdimen8\font{\vss\box0}\hss}%
     \or\hbox{\vrule width\wd2 height\fontdimen8\font depth0pt\hss}
     \or\setbox0=\hbox{a}
	\hbox to\ht0{\kern\dp0\vbox to\fontdimen8\font{\vss\box0}\hss}%
     \or\setbox0=\hbox{d}%
	\hbox to\ht0{\vrule width\fontdimen8\font height\ht2 depth\dp2\hss}
     \or\setbox0=\hbox{d}%
	\hbox to\ht0{\hss\vrule width\fontdimen8\font height\ht2 depth\dp2}
     \or\setbox0=\hbox{c}
	\hbox to\ht0{\vbox{\copy0\kern-\dp0}\hss}%
     \or\hbox{\vrule width\wd2 height\fontdimen8\font depth0pt\hss}
     \or\setbox0=\hbox{b}
	\hbox to\ht0{\kern\dp0\vbox{\copy0\kern-\dp0}\hss}%
     \or\box2 
     \fi
   }}


\newbox\brd_txt
\newbox\brd_box

\def\setborder[#1,#2,#3]#4
   {\def\brd_inbox
       {\hsize=#1
	\advance\hsize by-#2
	\advance\hsize by-#2
	\aftergroup\set_border
       }%
    \def\set_border
       {\setbox\brd_txt=\vbox
	 {\vskip#3
	  \box\brd_txt
	  \vskip#3
	 }%
	\setbox\brd_box=#4{\ht\brd_txt}{#1}%
	\vbox to\ht\brd_box
	 {\copy\brd_box \kern-\ht\brd_box
	  \vss
	  \hbox to\wd\brd_box{\hss\box\brd_txt\hss}
	  \vss
	 }%
       }%
    \afterassignment\brd_inbox
    \setbox\brd_txt=\vbox
   }


\catcode`\_=\borderelm 

\let\borderelm=\roundcorners 